\documentclass[11pt]{article}
\usepackage{amsmath,amssymb,verbatim}
\usepackage{times}

\usepackage[dvips]{graphicx}
\usepackage{amsfonts}
\usepackage{amsthm}
\usepackage{epsfig}
\usepackage{color}

\theoremstyle{plain}
\newtheorem{theorem}{Theorem}[section]
\newtheorem{lemma}[theorem]{Lemma}
\newtheorem{proposition}[theorem]{Proposition}
\newtheorem{corollary}[theorem]{Corollary}
\theoremstyle{definition}
\newtheorem{definition}[theorem]{Definition}

\newtheorem{example}[theorem]{Example}

\newtheorem{discussion}[theorem]{Discussion}

\newcommand{\todo}[1]{\vspace{5 mm}\par \noindent
\marginpar{\textsc{ToDo}}
\framebox{\begin{minipage}[c]{0.95 \textwidth}
\tt #1 \end{minipage}}\vspace{5 mm}\par}

\renewcommand{\todo}[1]{}

\newcommand{\idiot}[1]{\vspace{5 mm}\par \noindent
\framebox{\begin{minipage}[c]{0.95 \textwidth}
\tt #1 \end{minipage}}\vspace{5 mm}\par}

\renewcommand{\idiot}[1]{}

\newcommand{\xs}{x_1,\ldots,x_n}                
\newcommand{\Ms}{M_1,\ldots,M_q}                

\newcommand{\depth}{{\rm{depth}}\ }             
\newcommand{\dimn}{ {\rm{dim}} \ }
\newcommand{\be}{\beta}                          
\newcommand{\ddiv}{\mathrel{|}}
\newcommand{\ndiv}{\nmid}
\newcommand{\height}{ {\rm{height}} \ }         
\newcommand{\alx}[1]{{#1{}^{\vee}}}               
\newcommand{\st}{\mid}

\newcommand{\erase}[1]{}

\newcommand{\Il}{{\mathcal{I}}_{\mu}}
\newcommand{\Jl}{{\mathcal{J}}_{\mu}}

\newcommand{\El}{{\mathcal{E}}_{\mu}}
\newcommand{\J}{{\mathcal{J}}}
\newcommand{\E}{{\mathcal{E}}}
\newcommand{\pp}{{\mathfrak{p}}}             
\newcommand{\mm}{{\mathfrak{m}}}             
\newcommand{\set}{{\rm set}}

\newcommand{\ba}{ \,| \, }   
\newcommand{\FF}{\mathbb{F}}                
\newcommand{\reg}{{\rm reg}}

\def\Young#1{\vbox{\smallskip\offinterlineskip
    \halign{&\vbox{##}\kern-\Thickness\cr #1}}}

\newdimen\Squaresize \Squaresize=12pt
\newdimen\Thickness \Thickness=.3pt
\newdimen\Correction \Correction=7pt

\def\Vide#1{\hbox{
       \vbox to \Squaresize{\vss
          \hbox to \Squaresize{\hss#1 \hss}\vss}
    \hskip-\Correction}
   \kern-\Thickness}

\def\Carre#1{\hbox{\vrule width \Thickness
   \vbox to \Squaresize{\hrule height \Thickness\vss
      \hbox to \Squaresize{\hss#1\hss}
   \vss\hrule height\Thickness}
   \unskip\vrule width \Thickness}
   \kern-\Thickness}

\def\Box#1{\Carre{$\scriptstyle#1$}}

 \date{\today}

\author{Riccardo Biagioli\thanks{Dipartimento di Matematica, Universit\`a di Roma ``La Sapienza", biagioli@mat.uniroma1.it} \and
Sara Faridi\thanks{Department of Mathematics, Dalhousie University,
faridi@mathstat.dal.ca (research supported by NSERC)} \and Mercedes
Rosas\thanks{Department of Mathematics and Statistics, York
University, rosas@mathstat.yorku.ca}}

\title{Resolutions of De Concini-Procesi ideals indexed by hooks}

\begin{document}

\maketitle

\begin{abstract}
We find a minimal generating set for the De Concini-Procesi ideals
 indexed by hooks, and study their minimal free resolutions as well as
 their Hilbert series and regularity.
\end{abstract}


\section{Introduction}
In their study of the cohomology ring of the flag variety, De
Concini and Procesi~\cite{DP}, defined for any partition $\mu$ of
$n$, an ideal $\mathcal{I}_\mu$ of the polynomial ring
$R=\mathbb{Q}[\xs]$. In particular, they showed that the cohomology
ring of the variety of the flags fixed by a unipotent matrix of
shape $\mu$ may be presented as the graded quotient of the
polynomial ring $R$ by the ideal $\Il$. The space $R/\Il$ is
actually an interesting graded representation of the symmetric group
$S_n$, and it has been studied from different points of view by
several authors. Garsia and Procesi~\cite{GP}, studied its graded
character and showed that it can be expressed in terms of
Kostka-Foulkes polynomials. These polynomials appear in the
expansion of the classical Hall-Littlewood polynomials in the basis
of Schur functions \cite[Chapter III]{M2}, and were conjectured to
have positive integer coefficients. The result of Garsia and Procesi
mentioned above, gave an elegant proof of this positivity
conjecture.  N. Bergeron and Garsia \cite{BG} showed that as
symmetric group representations, the $R/\Il$ are isomorphic to
certain spaces of harmonic polynomials. Aval and N. Bergeron in
\cite{AB}, and Tanisaki in \cite{T} gave different sets of
generators for the ideal $\mathcal{I}_\mu$.  Finally, another
important feature of the $S_n$-modules $R/\mathcal{I}_\mu$ is that
they led Garsia and Haiman~\cite{GH} to the definition of the doubly
graded modules that appears in the famous $n!$ conjecture,  recently
solved by Haiman~\cite{H1}.

Despite the spaces $R/\Il$ having been extensively studied from the
point of view of representation theory and combinatorics, no commutative algebra investigation of these objects has been done so far. The goal of this paper is to begin that study. One of the strongest tools for finding
numerical informations about an ideal in a polynomial ring is finding
its minimal free resolution. The resolution in particular produces all
the numerical invariants that are described by the Hilbert function of
the ideal. Finding an exact description of the resolution for a
general ideal is usually a difficult task, there is a lot of
research and numerous conjectures on this problem. However, when the partition
$\mu$ indexing the De Concini-Procesi ideal $\Il$ is a hook, we are able to produce a minimal generating set for $\Il$, that we break into two parts. We show that one part forms an ideal with linear quotients, whose
resolutions are described by Herzog and
Takayama~\cite{HT}. We then show that the second part forms a regular sequence over the first part, and hence the resolution of this part is also well
understood. Below we describe this construction in detail, and compute
the Poincar\'e series associated to such an ideal (i.e. the generating
function encoding the ranks of the free modules appearing in a minimal
free resolution of the ideal). We also give a description of the
Hilbert series of $R/\Il$.

This paper is organized as follows. In Section~\ref{s:dp-ideals}, we
give the basic definitions of partitions and the language used in the
paper. We introduce De Concini-Procesi ideals, and compute a new
generating set for them in the case of hooks; we show later in
Section~\ref{ss:resolutions} that this generating set is
minimal. Section \ref{s:dp-res} contains a review of resolutions,
Cohen-Macaulay rings, and the other commutative algebra tools that we
use in the paper. In Section~\ref{ss:resolutions} we study the
resolutions of such ideals, and conclude with the formula of the
corresponding bigraded Poincar\'e series.  Finally, in
sections~\ref{s:regularity} and~\ref{s:hilbert}, we compute the
regularity and build the Hilbert series of the module $R/\Il$.

We hope that our exposition will appeal to readers not only in
commutative algebra, but also in combinatorics and invariant theory.
This is why, throughout the paper, we review the background material
we need from each field to make the concepts accessible to a wider
audience.

\noindent \emph{Acknowledgments:} All the test examples that supported
this research were run using the computer algebra program
Macaulay2~\cite{Macaulay2}. We would like to thanks Fran\c cois
Bergeron and Tony Geramita for constant support and for useful
discussions.


\section{De Concini-Procesi Ideals}\label{s:dp-ideals}

In this section, we introduce a family of ideals of the polynomial ring $R=k[\xs]$ indexed by partitions of $n$. These ideals were first introduced by De Concini-Procesi \cite{DP}, as ideals of the polynomial ring with rational coefficients. For our purpose $k$ may be an arbitrary field of characteristic $0$.
Let us start with some definitions and notation about partitions, that will be
used in the rest of this paper.

\subsection{Partitions}\label{s:partitions}

We let $\mathbb{P}=\{1,2,\ldots\}$, and $\mathbb{N}=\mathbb{P}\cup \{0\}$.
The cardinality of a set $S$ is denoted by $|S|$.
We define a {\em partition} of $n \in \mathbb{N}$ to be a finite sequence
$\mu=(\mu_1,\ldots,\mu_k) \in \mathbb{N}^k$, such that $\sum_{i=1}^k \mu_i =n$ and $\mu_1 \geq \ldots \geq \mu_k$. If $\mu$ is a partition of $n$ we write $\mu \vdash n$. The nonzero terms $\mu_i$ are called {\em parts} of $\mu$. The number
of parts of $\mu$ is called the {\em length} of $\mu$, denoted
by $\ell(\mu)$.

\begin{figure}[h]
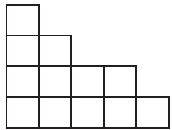

\begin{eqnarray*}
&\begin{matrix} \Young{ \Box{}\cr \Box{}& \Box{} \cr \Box{}& \Box{}&
               \Box{}& \Box{}\cr \Box{}& \Box{}& \Box{}& \Box{}&
               \Box{} \cr }
     \end{matrix}
 \end{eqnarray*}
\caption{The partition $\mu=(5,4,2,1)$}\label{Ytableau}
\end{figure}

The {\em Young diagram} of a partition  $(\mu_1,\ldots,\mu_k)\vdash
n$, is the diagram with $\mu_i$ squares in the $i^{\rm th}$-row. We
use the symbol $\mu$ for both a partition and its associated Young
diagram. For example, the diagram of $\mu=(5,4,2,1)$ is illustrated
in Figure \ref{Ytableau}.

For a partition $\mu=(\mu_1,\ldots,\mu_k)$
denote the {\em conjugate} partition $\mu^\prime:=(\mu_1^\prime,\ldots,
\mu_h^\prime)$, where for each $i\geq 1$, $\mu_i^\prime$ is the number of
parts of $\mu$ that are bigger than or equal to $i$.
The diagram of $\mu^\prime$ is obtained by flipping the diagram of $\mu $ across the diagonal.

Partitions of the form $\mu=(a)$ and $\mu=(1^b)=(\overbrace{1,\ldots,1}^b)$,  with $a,b \in \mathbb{P}$ are called {\em one-row} and  {\em one-column} partitions,
respectively. More generally, a partition is said to be a {\em hook}
if it is of the form $\mu=(a+1,1^b)$, with $a,b \in \mathbb{N}$.

\begin{figure}[h]
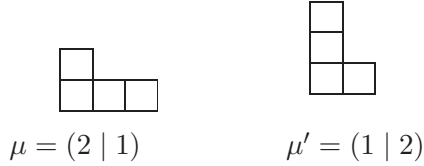

\begin{eqnarray*}
&\begin{matrix} \Young{ \; \cr \; \cr \; \cr\; \cr\; \cr \; \cr\; \cr
               \cr \; \cr \Box{}\cr \Box{} & \Box{} &\Box{}  \cr }
     \end{matrix}  \;\;\;\;\;\;\;\;\;\;\;\;\;\;\;\;\;\;\;
     \begin{matrix}
               \Young{ \Box{}\cr \Box{}\cr \Box{} & \Box{}
               \cr}
     \end{matrix}\\
    & \mu=(2 \mid 1) \;\;\;\;\;\;\;\;\;\;\;\;\;\;\;\;\;\;
    \mu^\prime=(1\mid 2)
 \end{eqnarray*}
\caption{Frobenius notation}\label{Frobenius-notation}
\end{figure}

Sometimes, it will be useful to denote hook partitions using a different notation. The hook $\mu=(a+1,1^b)$ in  {\em
Frobenius's notation} \cite[page 3]{M2} will be denoted by $\mu=( a \mid b )$. Note that its conjugate is $\mu^\prime= (b
\mid a)$. See Figure \ref{Frobenius-notation} for an example.

\subsection{De Concini-Procesi ideals}

From now on, we shall assume that a partition of $n$ has $n$ terms. So we will add enough zero terms to any partition until we have the right number of terms.
Let $\mu=(\mu_1,\ldots,\mu_n)$ be a partition of $n$, and
$\mu^\prime=(\mu^\prime_1\ldots,\mu^\prime_n)$ its  conjugate
partition.  For any $1\leq k \leq n$,
we define
\[\delta_k(\mu):=\mu_n^\prime + \mu_{n-1}^\prime + \ldots + \mu_{n-k+1}^\prime.\]

Recall that for any $1 \leq r \leq n$, the {\em elementary symmetric
polynomial} \cite{St2} is defined by
\[e_r(\xs):=\sum_{1\leq i_1<\ldots<i_r \leq n} x_{i_1}x_{i_2}\cdots x_{i_r}.\]
Given a subset $S \subseteq \{\xs\}$,  let $e_r(S)$ be the $r^{\rm
th}$ elementary symmetric polynomial in the variables in
$S$. Clearly, every $e_r(S)$ is a homogeneous polynomial in $R$ of
degree $r$.

We are now ready to introduce the ideals originally defined by De Concini and Procesi \cite{DP}. We use a different and simpler set of generators which was defined by Tanisaki \cite{T}.

\begin{definition}[De Concini-Procesi ideal] We let $\mathcal{C}_\mu$ denote the
collection of partial elementary symmetric polynomials
\begin{equation}\label{generators}
\mathcal{C}_\mu=\{e_r(S) \mid S\subseteq \{\xs\}, \; |S|=k \geq 1, \; k \geq r > k -\delta_k(\mu)
 \}.
\end{equation}
The {\em De Concini-Procesi ideal} $\mathcal{I}_\mu$ is  the
homogeneous ideal generated by the elements of $\mathcal{C}_\mu$, in
symbols,
\begin{equation*}
\mathcal{I}_\mu:=(\mathcal{C}_\mu).
\end{equation*}
\end{definition}
Note that $\delta_n(\mu)=n$, for any partition $\mu$ of $n$. Hence when we set $k=n$ in (\ref{generators}), we obtain that $\mathcal{I}_\mu$ contains the ideal generated by the elementary symmetric polynomials in all the
variables $x_1,\ldots,x_n$.

\begin{example}\label{ex1} Let $\mu=(3,1,0,0)\vdash 4$ and $\mu^\prime=(2,1,1,0)$ be the partitions appearing in Figure \ref{Frobenius-notation}. Then   $(\delta_1(\mu),\ldots,\delta_4(\mu))=(0,1,2,4)$. Hence
\[(1-\delta_1(\mu),\ldots,4-\delta_4(\mu))=(1,1,1,0),\]
and the collection $\mathcal{C}_\mu$ consists of the following
elements.  For $k=1$ there is no admissible $e_r(S)$. For $k=2$ we
get the set of monomials:
\[x_1 x_2, \; x_1 x_3, \; x_1 x_4, \; x_2x_3, \; x_2x_4, \; x_3x_4.\]
For $k=3$:
\[x_1 x_2 + x_1 x_3 + x_2 x_3, \; x_1x_2 + x_1x_4 + x_2 x_4, \; x_1 x_3 + x_1 x_4 + x_3 x_4, \; x_2 x_3 + x_2 x_4 + x_3 x_4,  \]
\[ x_1 x_2 x_3,  \; x_1 x_2 x_4, \; x_1 x_3 x_4, \; x_2 x_3 x_4.\]
Finally for $k=4$, as already noted, we get the complete set of the elementary symmetric
functions $e_r(x_1,x_2,x_3,x_4)$, for $1 \leq r \leq 4$.
\end{example}

When the indexing partition $\mu$ is a hook, the ideal
$\mathcal{I}_\mu$ can be split in two parts. We have the following
result.

\begin{proposition}[Generators of De Concini-Procesi ideals indexed by hooks]\label{p:generators}
Let $\mu=(a \mid b)\vdash n $ be a hook.
Then the De Concini-Procesi ideal associated to $\mu$ in the polynomial ring $k[\xs]$ is
$$\Il=\Jl + \El,$$ where
\begin{eqnarray}\label{e:Jdef}
\Jl=(x_{i_1}\cdots x_{i_{b+1}} \mid 1\leq i_1 < \ldots< i_{b+1} \leq n )
\end{eqnarray}
is the ideal generated by all
square-free monomials in $\xs$ of degree $b+1$, and
\begin{eqnarray}\label{e:Edef}
\El=( e_i(x_1,\ldots,x_n) \mid 1\leq i \leq b)
\end{eqnarray}
is the ideal generated by all elementary symmetric polynomials of degree $\leq b$
in the variables $\xs$.
\end{proposition}

\begin{proof}  The partition $\mu=(a \mid b)$ is of size $n=a+b+1$. We can write $$\mu^\prime=(b\mid a)=(b+1,\underbrace{1,\ldots,1}_a,\underbrace{0,\ldots,0}_b).$$
Then we have
$$
(\delta_1(\mu), \delta_2(\mu), \ldots, \delta_n(\mu))=(\underbrace{0, \ldots, 0}_b, 1, 2, \ldots, a, n),$$ and so
$$(1-\delta_1(\mu), 2-\delta_2(\mu), \ldots, n-\delta_n(\mu))=(1,2,3, \ldots, b, \underbrace{b,\ldots,b}_a,0).
$$
The definition of $\mathcal{C}_\mu$ in (\ref{generators}) implies that no $k$, with $1\leq k \leq b$, contributes  a generator to the ideal $\Il$.

The first index making a  nontrivial contribution to the set $\mathcal{C}_\mu$ is
$k=b+1$, which adds to $\mathcal{C}_\mu$ all $e_{b+1}(S)$, with $|S|=b+1$, or in other words all the square-free monomials of degree $b+1$ in the variables $\xs$. We denote by $\Jl$ the ideal generated by these square-free monomials.

Now all the indices $k$, with $b+2 \leq k \leq n-1$ add to $\mathcal{C}_\mu$ elements of the form $e_r(S)$, with $k\geq r\geq b+1$, and $|S|=k$. Each such $e_r(S)$ is a homogeneous polynomial of degree $r$, which we can write as the sum of square-free monomials of degree $r$. Since $r \geq b+1$, and all square-free monomials of degree $b+1$ or more are already in $\Il$, such $e_r(S)$ do not contribute any new generators to $\Il$.

Finally, for $k=n$ we obtain all the elementary symmetric polynomials in all the variables. For the same reasons as above, the only new contributions are
\begin{equation*}
e_1(x_1, \ldots, x_n), e_2(x_1, \ldots, x_n),\ldots, e_b(x_1, \ldots,
x_n).
\end{equation*}
We denote the ideal generated by these elementary symmetric polynomials by $\El$. We conclude that $\Il=\Jl + \El$
\end{proof}

\begin{example} Let $\mu=(2 \mid 1)\vdash 4$. It follows from the computations in Example \ref{ex1}, that the ideal $\Il$ splits into two parts
$$\Il=(x_1 x_2 ,x_1 x_3,x_1 x_4 ,x_2x_3,x_2x_4,x_3x_4) +
(x_1+x_2+x_3+x_4).$$  The first part is generated by
all monomials of degree 2 in the variables $x_1,x_2,x_3,x_4$, and
the second is generated by $e_1(x_1,x_2,x_3,x_4)$, the elementary symmetric polynomial of degree 1.
\end{example}


\section{Commutative algebra tools}\label{s:dp-res}

\todo{(Sara)
\begin{itemize}

\item Are the results  characteristic free? Probably, but since the
original result of DP was in char 0, we'll leave it there for now.

\item Most results are for local rings, but *-local is
enough. Herzog's Lectures~\cite{He}, for example, have been done for
graded $k$-algebras. Also see~\cite[Section~1.5]{BH} to see how
everything corresponds between local and *-local.

\item (TO DO) Find Hilbert polynomial!  [tried-- this is hard! But not
impossible.... There are some neat relationships.

\item (TO DO) In the final reading, make sure all the arguments
between resolutions of $I$ and $R/I$ match.

       \end{itemize}
}

Let $R=k[\xs]$ be a polynomial ring over a
field $k$ of characteristic 0, with the standard grading $\deg x_i=1$,
for all $i$. Let $\mm=(\xs)$ be the (irrelevant) homogeneous maximal
ideal of $R$. We are usually interested in the quotient $S=R/I$ where
$I$ is an ideal of $R$ generated by homogeneous polynomials. In this
situation, $S$ inherits the grading and the irrelevant maximal ideal
from $R$ via the quotient map. Much of what we discuss below will be
in this context, but applies more generally to local rings.

\subsection{Resolutions}

Resolutions provide us with an effective method to study a finitely
generated module $M$ via a sequence of free modules mapping to
it. Among the many applications, the ranks of these free modules, also
known as ``Betti numbers'', are numerical invariants of $M$ that make
it possible to compute the Hilbert function of $M$ directly.  There is
a large amount of literature focusing on different aspects of
resolutions, studying them using homological, geometrical, or
combinatorial tools. We refer the interested reader
to~\cite{E2},~\cite{He} or~\cite{BH} to learn more about current
research in this field. Eisenbud's book \cite{E2} in particular
contains a beautiful exposition on the history of the subject.

\begin{definition}[Minimal free resolution]\label{d:resolution}
A \emph{free resolution} of $R/I$ is an exact complex $\FF$
 $$0 \longrightarrow \cdots
                \stackrel{\delta_{i+1}}{\longrightarrow}
                F_i
                \stackrel{\delta_i}{\longrightarrow}
                F_{i-1}
                \stackrel{\delta_{i-1}}{\longrightarrow} \cdots
                \stackrel{\delta_{2}}{\longrightarrow}
                F_1
                \stackrel{\delta_{1}}{\longrightarrow}
                R
                \stackrel{\delta_{0}}{\longrightarrow} R/I
                \longrightarrow 0.
             $$
of free $R$-modules $F_i$ ($F_0=R$). The resolution is \emph{minimal} if  $\delta_i(F_i) \subseteq \mm F_{i-1}$ for $i>0$.
\end{definition}

It is worth noting that the difference between the resolution of the
ideal $I$ (as an $R$-module), and the resolution of the quotient $R/I$
is just the one free module $F_0$: given the
above resolution for $R/I$, the resolution for $I$ is the following:
$$0 \longrightarrow \cdots
                \stackrel{\delta_{i+1}}{\longrightarrow}
                F_i
                \stackrel{\delta_i}{\longrightarrow}
                F_{i-1}
                \stackrel{\delta_{i-1}}{\longrightarrow} \cdots
                \stackrel{\delta_{2}}{\longrightarrow}
                F_1
               \stackrel{\delta_{1}}{\longrightarrow} I
                \longrightarrow 0.
             $$
In this paper, we will always be considering the resolution of $R/I$.

It follows from the Hilbert Syzygy Theorem \cite[Theorem 1.1]{E2} that
the length of a minimal free resolution of $R/I$ is finite; i.e. $F_i=0$
for $i>n$, where $n$ is the number of variables in $R$ (the resolution
could stop even earlier, there are formulas to compute the length of a
resolution).  A minimal free resolution of $R/I$ is unique up to
isomorphism~\cite[Theorem 20.2]{E1}.


If each $F_i$ is a free module of rank $\be_i$, the resolution of $R/I$ is
 \begin{eqnarray}\label{e:first-resolution}0 \longrightarrow
                R^{\beta_m}
                \stackrel{\delta_{m}}{\longrightarrow}
                R^{\beta_{m-1}}
                \stackrel{\delta_{m-1}}{\longrightarrow}
                \cdots
                \stackrel{\delta_{2}}{\longrightarrow}
                R^{\beta_1}
                \stackrel{\delta_{1}}{\longrightarrow} R
                \stackrel{\delta_{0}}{\longrightarrow} R/I
                \longrightarrow 0.
 \end{eqnarray}
The $\beta_i$ are called the \emph{Betti numbers} of $R/I$; these are
independent of which minimal resolution one considers.

In the case where $I$ is a homogeneous ideal, and therefore $R/I$ is
graded, we define the \emph{graded Betti numbers} of $R/I$.  This is done by making the maps $\delta_i$ homogeneous, so that they take a degree $j$ element of $F_i$ to a degree $j$ element of $F_{i-1}$. To serve this purpose the degree of each generator of $F_i$ is adjusted.  So we can write the free module $F_i=R^{\beta_i}$
as
$$R^{\beta_i}=\bigoplus_{j}R(-j)^{\beta_{i,j}}$$ where for a given
             integer $a$, $R(a)$ is the same as $R$ but with a new
             grading:
             $$R(a)_d=R_{a+d}.$$ So the resolution shown in
             (\ref{e:first-resolution}) becomes

\begin{eqnarray}\label{e:second-resolution}0 \longrightarrow
               \bigoplus_{j}R(-j)^{\beta_{m,j}}
                \stackrel{\delta_m}{\longrightarrow}
               \bigoplus_{j}R(-j)^{\beta_{m-1,j}}
                \stackrel{\delta_{m-1}}{\longrightarrow} \cdots
                \stackrel{\delta_{2}}{\longrightarrow}
               \bigoplus_{j}R(-j)^{\beta_{1,j}}
                \stackrel{\delta_{1}}{\longrightarrow} R
                 \stackrel{\delta_{0}}{\longrightarrow} R/I
                 \longrightarrow 0.
 \end{eqnarray}
This is called the \emph{graded minimal free resolution} of $R/I$,
and the $\beta_{i,j}$ are the \emph{graded Betti numbers} of
$R/I$. Clearly, $\displaystyle \sum_j \beta_{i,j}=\beta_i$.

\begin{definition}[Bigraded Poincar\'e series]\label{d:poincare}
The \emph{bigraded Poincar\'e series} of an ideal $I$ is the
generating function for the graded Betti numbers of $I$:
$$P_{R/I}(q,t)=\sum_{i,j}\beta_{i,j}q^it^j.$$

\end{definition}

\begin{definition}[Linear resolution] The graded resolution described
 in (\ref{e:second-resolution}) is a \emph{linear resolution}, if for
some $u$, $\beta_{i,j}=0$ unless $j=u+i-1$. In other words, $R/I$ has
a linear resolution if for some $u$, it has a graded minimal free
resolution of the form

$$              0 \longrightarrow
               R(-(u+m-1))^{\beta_{m, u+m-1}}
                \stackrel{\delta_{m}}{\longrightarrow}
                R(-(u+m-2))^{\beta_{m-1,u+m-2}}
                \stackrel{\delta_{m-1}}{\longrightarrow}
$$
$$                \cdots
                \stackrel{\delta_{3}}{\longrightarrow}
                R(-(u+1))^{\beta_{2,u+1}}
               \stackrel{\delta_{2}}{\longrightarrow}
                R(-(u))^{\beta_{1,u}}
                \stackrel{\delta_{1}}{\longrightarrow} R
                \stackrel{\delta_{0}}{\longrightarrow} R/I
                \longrightarrow 0.$$
In this case, all the generators of the ideal $I$ have degree equal to
$u$.
\end{definition}

\begin{discussion}[Resolutions using mapping cones]\label{d:mapping-cones}
The mapping cone technique provides a way to build a free resolution
of an ideal by adding generators one at a time.  A resolution
obtained using mapping cones is not in general minimal. However, we
will be focusing only on the special case of multiplication by a
nonzerodivisor, in which case we obtain a minimal free resolution. For
a more general or detailed description, see~\cite{Sc},~\cite{HT},
or~\cite{E2}.

Suppose that $I$ is an ideal in the polynomial ring
$R$, and $e \in \mm$ is a nonzerodivisor in $R/I$ (i.e. $e$ is a
regular element mod $I$; see Definition~\ref{d:depth}). The goal is
to build a minimal free resolution of $R/I+(e)$ starting from a
minimal free resolution of $R/I$. Consider the short exact
sequence
$$0 \longrightarrow R/I:(e) \stackrel{.e}{\longrightarrow} R/I
\longrightarrow R/I+(e) \longrightarrow 0$$ where $I:(e)$ is the
quotient ideal consisting of all elements $x \in R$ such that $xe \in
I$.  Since $e$ is a nonzerodivisor in $R/I$, we have $I:(e)=I$, and so
our short exact sequence turns into
$$0 \longrightarrow R/I \stackrel{.e}{\longrightarrow} R/I
\longrightarrow R/I+(e) \longrightarrow 0.$$

Suppose we have a minimal free resolution of $R/I$
            \begin{eqnarray}\label{e:resofR/I}
            0 \longrightarrow \cdots
                \stackrel{\delta_{i+1}}{\longrightarrow}
                A_i
                \stackrel{\delta_i}{\longrightarrow}
                A_{i-1}
                \stackrel{\delta_{i-1}}{\longrightarrow} \cdots
                \stackrel{\delta_{2}}{\longrightarrow}
                A_1
                \stackrel{\delta_{1}}{\longrightarrow} R
                          \stackrel{\delta_{0}}{\longrightarrow} R/I
                \longrightarrow 0.
             \end{eqnarray}
Then we can obtain the following minimal free resolution of $R/I+(e)$
            \begin{eqnarray}\label{e:mappingcone}
              0 \longrightarrow \cdots
                \stackrel{d_{i+1}}{\longrightarrow}
                F_i
                \stackrel{d_i}{\longrightarrow}
                F_{i-1}
                \stackrel{d_{i-1}}{\longrightarrow} \cdots
                \stackrel{d_{2}}{\longrightarrow}
                F_1
                \stackrel{d_{1}}{\longrightarrow} R
                          \stackrel{d_{0}}{\longrightarrow} R/I +(e)
                \longrightarrow 0
             \end{eqnarray}
where for each $i>0$, as a free \erase{(not necessarily graded)} $R$-module
\begin{center}$F_i=A_i\oplus A_{i-1}$ and $d_i(x,y)=(ey+\delta_i(x),
-\delta_{i-1}(y))$.\end{center}
The resolution is minimal because for $(x,y) \in A_i\oplus
A_{i-1}$ and $e \in \mm$, we have $$ey \in \mm A_{i-1}, \ \delta_i(x) \in
\mm A_{i-1},\ \delta_{i-1}(y) \in \mm A_{i-2} \Longrightarrow
d_i(x,y)\in \mm F_{i-1}.$$

We now focus on the grading of each $F_i$.  Suppose that the element
$e \in R$ is homogeneous of degree $m$, and for each $i$, each of the free modules $A_i$ in (\ref{e:resofR/I}) are of the form
$$A_i=\bigoplus_{j}R(-j)^{\beta_{i,j}}$$ where the $\beta_{i,j}$ are
the graded Betti numbers. We would like to compute the graded Betti
numbers of $R/I+(e)$.

\begin{lemma}\label{l:mapping-cone} Let $I$ be an ideal of
 the polynomial ring $R$, and $e \in \mm$ a homogeneous element of
degree $m$ which is a  nonzerodivisor in $R/I$. Consider the minimal
free resolutions (\ref{e:resofR/I}) of $R/I$, and (\ref{e:mappingcone})
of $R/I+(e)$ obtained by mapping cones. For each $i>0$ we have
$$F_i=\bigoplus_{j}R(-j)^{\beta_{i,j}}\ \oplus \
\bigoplus_{j}R(-j-m)^{\beta_{i-1,j}}.$$
\end{lemma}

           \begin{proof} We prove this by induction on $i$. In the case
             where $i=1$, we have the homogeneous map
             $$d_1: A_1 \oplus R \longrightarrow R$$ where
             $d_1(x,y)=ey+\delta_1(x)$. In particular, if $x \in A_1$
             is a homogeneous element of degree $t$, then
             $d_1(x,0)=\delta_1(x)$ is also a degree $t$ homogeneous
             element of $R$. If $y \in R$ is a homogeneous element of
             degree $t$, then $d_1(0,y)=ey$ has degree $t+m$. In
             order to make $d_1$ a homogeneous (degree 0) map, we
             shift the grading of the component $R$ of $F_i$ by $m$, so that
             $$F_1=\bigoplus_{j}R(-j)^{\beta_{1,j}}\oplus R(-m).$$

             Suppose our claim holds for all indices less than $i$,
             and we have the homogeneous map $$ d_i: F_i=A_i\oplus
             A_{i-1} \longrightarrow
             F_{i-1}=\bigoplus_{j}R(-j)^{\beta_{i-1,j}}\ \oplus \
             \bigoplus_{j}R(-j-m)^{\beta_{i-2,j}}.$$ We use the same
             argument as we did in the case of $i=1$. If $x \in A_i$
             is a homogeneous element of degree $t$, then
             $d_i(x,0)=\delta_i(x)$ is also a degree $t$ homogeneous
             element of $A_{i-1}$. If $y \in A_{i-1}$ is a homogeneous
             element of degree $t$, then
             $d_i(0,y)=(ey,-\delta_{i-1}(y))$ has to be a homogeneous
             element of $F_{i-1}$ of degree $t$. By definition, this
             is already true for the component $\delta_{i-1}(y)$, but
             $ey$ has degree $m+t$. So in order to make $d_i$ a
             homogeneous (degree 0) map, we have to shift the grading
             of each component of $F_i$ that comes from $A_{i-1}$ by
             $m$, so that
             $$F_i=\bigoplus_{j}R(-j)^{\beta_{i,j}}\ \oplus \
             \bigoplus_{j}R(-j-m)^{\beta_{i-1,j}}.$$

       \end{proof}

\end{discussion}

\begin{corollary}\label{c:mapping-cone} Let $I$ be an ideal of
 the polynomial ring $R$ and $e \in \mm$ be a homogeneous element of
degree $m$ which is a nonzerodivisor in $R/I$. Then
$$P_{R/I+(e)}(q,t)= (1+qt^m)P_{R/I}(q,t).$$
\end{corollary}

         \begin{proof} By Lemma~\ref{l:mapping-cone}, if for a fixed $i$,
          $\displaystyle A_i=\bigoplus_{j=0}^{b_i}
            R(-j)^{\beta_{i,j}}$ then $$\displaystyle
            F_i=\bigoplus_{j=0}^{b_i}R(-j)^{\beta_{i,j}}\ \oplus \
            \bigoplus_{j=0}^{b_{i-1}}R(-j-m)^{\beta_{i-1,j}}.$$

           So we have
           $$\begin{array}{ll}
            P_{R/I+(e)}(q,t)& \displaystyle =1 + \sum_{i\geq 1} \left
            (\sum_{j=0}^{b_i} \beta_{i,j}t^j +
            \sum_{j=0}^{b_{i-1}}\beta_{i-1,j}t^{j+m} \right ) q^i\\

            & \displaystyle = \sum_{i\geq 0} \sum_{j=0}^{b_i}
            \beta_{i,j}t^jq^i +  t^m \sum_{i \geq 0}
            \sum_{j=0}^{b_i}\beta_{i,j}t^j q^{i+1}\\

            & \displaystyle = (1+qt^m)\sum_{i\geq 0} \sum_{j=0}^{b_i}
            \beta_{i,j}t^jq^i \\

            &\\

            & \displaystyle = (1+qt^m)P_{R/I}(q,t).

             \end{array}$$

     \end{proof}


\subsection{Krull dimension, Cohen-Macaulay rings, minimal primes}
A minimal prime
ideal (with respect to inclusion) containing $I$ is called a
\emph{minimal prime} of $I$.
Given any ideal $I$ of $R$, the \emph{Krull dimension} or
\emph{dimension} of the quotient ring $R/I$ is equal to the length $r$
of the maximal chain of prime ideals containing $I$
 $$\pp_0 \subset \pp_1 \subset \cdots \subset \pp_r$$
(here, $\pp_0$ is a minimal prime of $I$).
The \emph{height} of a prime ideal $\pp$ is the maximal length of a
chain of prime ideals
$$ \pp=\pp_0 \supset \pp_1 \supset \cdots \supset \pp_r,$$ and the height of
a general ideal $I$ is the smallest height of its minimal primes.

\begin{example} If $I=(xy,xz) \subset R=k[x,y,z]$, then
the minimal primes of $I$ are $(x)$ and $(y, z)$. In this case
$\height I=1$ (as $(x) \supset (0)$ is a maximal chain) and $\dimn
R/I= 2$. \end{example}

\begin{definition}[Regular sequence, depth, Cohen-Macaulay]\label{d:depth}
Let $R$ is a polynomial ring with standard grading, and $I$  a homogeneous ideal of $R$. Consider  $S=R/I$  with homogeneous maximal ideal $\mm$.
A sequence $y_1,\ldots,y_m$ of elements in $\mm$ is a \emph{regular sequence} of $S$ if

\begin{tabular}{ll}
(i) &$(y_1, \ldots, y_m)S\neq S$,\\
(ii)& $y_1$ is a nonzerodivisor in $S$, \\
(iii)& $y_i$ is a nonzerodivisor in $S/(y_1,\ldots,y_{i-1})$.
\end{tabular}

The length of a maximal regular sequence in $S$ is called the
\emph{depth} of $S$.

In general, the depth of $S$ is less than or equal to the
dimension of $S$, but in the case equality is obtained, i.e. $\depth(S)=\dimn(S)$, the ring is
\emph{Cohen-Macaulay}.
\end{definition}

For more on dimension theory  and on the theory of Cohen-Macaulay rings, see~\cite{E1}, Appendix A of~\cite{BH},
or~\cite{V}.

We will need the definition of the
dual of a square-free monomial ideal.  This is the same as Alexander
dual, but we state the (equivalent) definition in a slightly different
language (see~\cite{F} for more).
Recall that a \emph{monomial ideal} is an ideal generated by monomials, and a
\emph{square-free} monomial ideal is an ideal generated by square-free
monomials in the variables $\xs$.
\begin{definition}[dual of an ideal] Let $I$ be a square-free
monomial ideal in the polynomial ring $k[\xs]$. Then $\alx{I}$ is a
square-free monomial ideal, where each generator of $\alx{I}$ is the
product of the variables appearing in the generating set of a minimal
prime of $I$.
\end{definition}

Note that if $I$ is a monomial ideal, its minimal primes are generated
by single variables.

\begin{example} If $I=(xy,xz,yzw) \subset k[x,y,z,w]$, then
the minimal primes of $I$ are $(x, y)$, $(x, z)$, $(x, w)$ and $(y,
z)$. So $\alx{I}=(xy,xz,xw,yz)$. \end{example}

Recall that if $I$ and $J$ are two ideals of $R$, their
\emph{quotient} is the ideal defined as $$I:J=\{x\in R\st xJ
\subseteq I\}.$$

\begin{definition}[linear quotients]\label{d:linear-quotients} If
$I \subset k[\xs]$ is a monomial ideal and $G(I)$ is its unique
minimal set of monomial generators, then $I$ is said to have
\emph{linear quotients} if there is an ordering $\Ms$ on the elements
of $G(I)$ such that for every $i=2,\ldots,q$, the quotient ideal
$$(M_1,\ldots,M_{i-1}):M_i$$ is generated by a subset of the variables
$\xs$.
\end{definition}

\begin{lemma}\label{l:I-lemma} Let $I$ be an ideal in the polynomial
ring $R=k[\xs]$ generated by all square-free monomials of a fixed
degree $m$. Then

\begin{enumerate}

\item $I$ has linear quotients;

\item $R/I$ has a linear resolution;

\item $R/I$ is Cohen-Macaulay.

\end{enumerate}
\end{lemma}

      \begin{proof} We can order the generating monomials of
        $I$ lexicographically as $\Ms$. Take such a monomial
        $M_i=x_{j_1}\ldots x_{j_{m}}$, written so that $j_1 < j_2 <
        \ldots < j_{m}$.  Since $(M_1,\ldots,M_{i-1})$ is a monomial
        ideal, and $M_i$ is also a monomial, the quotient ideal
        $(M_1,\ldots,M_{i-1}):M_i$ is generated by monomials.
        Observe that

       \begin{enumerate}

       \item If $s< j_t$ for some $j_t \in \{j_1,\ldots, j_{m}\}$ and
     $s \not\in \{j_1,\ldots, j_{m}\}$, then $x_s \in
     (M_1,\ldots,M_{i-1}):M_i$.

       This is because the monomial $\frac{x_sM_i}{x_{j_t}}$ is a
       degree $m$ monomial that is lexicographically smaller than
       $M_i$, that is, $\frac{x_sM_i}{x_{j_t}} \in
       \{M_1,\ldots,M_{i-1}\}$.

      \item If $u$ is a monomial in $(M_1,\ldots,M_{i-1}):M_i$, then
       $M_l \ddiv uM_i$ for some $l <i$. Since $M_l <_{lex} M_i$, there
       exists $x_s$, such that $x_s\ddiv M_l$, $x_s \ndiv M_i$ and $s <j_t$
       for some $j_t \in \{j_1, \ldots, j_{m}\}$.

      \end{enumerate}

 It follows that $x_s \ddiv u$, and $(M_1,\ldots,M_{i-1}):M_i$ is
      generated by the set of variables $x_s$, with $s <j_m$ and $s \notin \{j_1,\ldots,j_m\}$ as described in part 1.
      This proves that $I$ has linear quotients.

      Since the generators of $I$ all have the same degree, and since
      $I$ has linear quotients, it also follows that $R/I$ has linear
      resolution (Lemma~5.2 of~\cite{F}).

      Now we focus on the structure of $\alx{I}$. Since every
      generator of $I$ has exactly $m$ variables, each such generator
      misses exactly $n-m$ variables from the set  $\{\xs\}$. So if
      $A$ is any $n-m+1$-subset of $\{\xs\}$, $A$ must contain at least
      one variable from each of the $M_i$. Also no proper subset of
      $A$ will have this property (i.e. $A$ is the minimal set with
      such a property). So $A$ is a generating set for a minimal prime
      of $I$. Since all minimal primes of $I$ are generated by subsets
      of $\{\xs\}$, it follows that $\alx{I}$ is generated by all
      square-free monomials of degree $n-m+1$.

      So we have shown that $\alx{I}$ satisfies the hypotheses of our
      lemma, and hence it satisfies parts 1 and 2. In particular,
      $R/\alx{I}$ has linear resolution, and so by  Theorem~3 of~\cite{ER}, equivalently,  $R/I$ is Cohen-Macaulay.

     \end{proof}

\section{Resolutions of De Concini-Procesi ideals of
hooks}\label{ss:resolutions}

In this section we study the minimal free resolutions of
the De Concini-Procesi ideal $\Il$ of a hook $\mu=(a\mid b)$. We have seen that $\Il$ is
the sum of two ideals $$\Il=\Jl+\El$$
where $\Jl$ is generated by monomials, and  $\El$  is generated by elementary symmetric functions. Below we show
how we can recover the resolution of $\Il$ using the resolutions
of each one of the summands.

Since  $\Jl$ is generated by all square-free monomials of
$R=k[\xs]$ that have degree $b+1$,  by Lemma~\ref{l:I-lemma},
$\Jl$ is a Cohen-Macaulay ideal with linear resolutions and linear
quotients. On the other hand, it is easy to see that all the minimal primes of
$\Jl$ have uniform height $n-b$. This is because every generator
of $\Jl$ is a product of exactly $b+1$ variables in the set $\{\xs\}$,
and so a minimal subset of $\{\xs\}$ that shares at least one variable
with each one of these generators must have $n-b$ elements. Such
an ideal will have height equal to $n-b$, and so it follows that
$\dimn R/\Jl=b$.

We have thus shown that

\begin{corollary}\label{c:J-CM} For a hook $\mu=(a \ba b)$,
the ideal $\Jl$ of $R$ has linear quotients, linear resolution, and
$R/\Jl$ is Cohen-Macaulay of (Krull) dimension $b$.
\end{corollary}

Next, we focus on the ideal $\El$, which is generated by the first $b$
elementary symmetric functions.

\begin{proposition}\label{p:e-regular}  For a hook $\mu=(a \ba b)$, the
 set of generators $$e_1(\xs), \ldots, e_b(\xs)$$ of $\El$ form a
regular sequence over the quotient ring $R/\Jl$. \end{proposition}

       \begin{proof}  Let $S=R/\Jl$. We know by Corollary~\ref{c:J-CM}
         that $S$ is a Cohen-Macaulay ring, and $\dimn S=b$. To show
         that $e_1(\xs), \ldots, e_b(\xs)$ forms a regular sequence in
         $S$, by Theorem 2.1.2 of~\cite{BH}, it is enough to show
         that $\dimn S/\El =0$.  We prove this by induction. Let
         $\underline{\mu}=((a-1) \ba b)$ and $\overline{\mu}=(a \ba (b-1))$ be two
         hooks consisting of $n-1$ squares. Notice that
         \begin{enumerate}
           \item $\Jl= x_n\J_{\overline{\mu}} + \J_{\underline{\mu}}$.

              To see this, split the generating set of $\Jl$ into
              two sets: $\underline{G}$ consists of all those generators that do
              not contain the variable $x_n$, and $\overline{G}$ is the
              rest. So $\underline{G}$ consists of all square-free monomials of
              degree $b+1$ in the variables $x_1, \ldots, x_{n-1}$,
              which is by definition the generating set of
              $\J_{\underline{\mu}}$.

              Similarly, if we factor out the variable
              $x_n$ from each monomial in $\overline{G}$, we will be left with
              all square-free monomials of degree $b$ in the variables
              $x_1, \ldots,x_{n-1}$, which by definition generate the
              ideal $\J_{\overline{\mu}}$.

          \item Since every term in $e_i(\xs)$ is a square-free
           monomial of degree $i$, we can partition all such monomials
           into those that contain $x_n$ and those that don't. It is
           then easy to see that for every $i$, $$e_i(\xs)=e_i(x_1,\ldots,x_{n-1}) +x_ne_{i-1}(x_1,\ldots,x_{n-1}).$$

     \end{enumerate}

         It follows that $$\frac{S}{\El +(x_n)}=\frac{k[\xs]}{\Jl +
         \El + (x_n)}\cong \frac{k[x_1, \ldots, x_{n-1}]}{\J_{\underline{\mu}} +
         \E_{\underline{\mu}}}$$ and so by the induction hypothesis,
         \begin{eqnarray} \label{e:dimension}\dimn \frac{S}{\El
         +(x_n)} =0.  \end{eqnarray} It follows that $\dimn S/\El =0$
         or 1.

     Now suppose that $\dimn S/\El= \dimn R/(\Jl
         +\El)=1$.

        \idiot{(Sara) We can assume that all rings are localized
         at the maximal ideal $\mm=(\xs)$. Since $R$ is a $^*$local
         ring, $R$ being Cohen-Macaulay is equivalent to $R_\mm$ being
         Cohen-Macaulay~\cite[Exercise 2.1.27(c)]{BH}.} So it follows
         that there is a prime ideal $\pp$ of $R$ such that
         $$\Jl+\El \subseteq \pp \subset \mm.$$

         Since $\pp$ is a prime ideal and every monomial generator of
          $\Jl$ belongs to $\pp$, at least one variable of $R$ has to
          be in $\pp$; say, $x_n \in \pp$ (the equality
          (\ref{e:dimension}) that we shall use holds if one replaces
          $x_n$ with any other variable in $R$).  But then $$\Jl + \El
          + (x_n) \subseteq \pp \subset \mm$$ but this contradicts the
          fact that $$\dimn \frac{k[\xs]}{\Jl + \El + (x_n)}=0.$$

        \idiot{(Sara)
        \begin{lemma} If $\pp$ is a prime ideal in $k[\xs]$, and
        the monomial $x_1\ldots x_r \in \pp$, then for some $i$, $x_i
        \in \pp$.\end{lemma}

        \begin{proof} Since    $x_1\ldots x_r \in \pp$,
        it follows that $x_1 \in \pp$ or $x_2\ldots x_r \in \pp$.  If
         $x_1 \in \pp$, we are done. Otherwise, we continue
         inductively until we have $x_{r-1}x_r \in \pp$, which implies
         that $x_{r-1} \in \pp$ or $x_r \in \pp$.\end{proof}}

        \end{proof}

 We are now ready to state our central claim.

\begin{theorem}[Main theorem]\label{mainthm} Let $\mu=(a \ba b)$ be a hook.
Then the bigraded Poincar\'e  series for the ideal $\Il$ is the
following
\begin{equation} \label{e:mainthm}
P_{R/\Il}(q,t)=\prod_{k=1}^{b} (1+qt^k) \,\cdot\, \Big( 1 + q
t^{b+1} \sum_{i=0}^a { b+i \choose b} (1+q\,t)^{i} \Big).
\end{equation}
\end{theorem}

       \begin{proof}  As usual, let
     $\Il=\Jl+\El$.

         \begin{itemize}

         \item[{\em Step 1.}] The ideal $\Jl$ has linear quotients
         (Corollary~\ref{c:J-CM}).  It follows from Corollary~1.6
         of~\cite{HT} that, if $G(\Jl)$ indicates the generating set
         for $\Jl$, the bigraded Poincar\'e series of $\Jl$ is the
         following:
          \begin{eqnarray}\label{e:poincare}
          P_{R/\Jl}(q,t)= 1 + \sum_{M \in G(\Jl)}
          (1+qt)^{|\set(M)|}qt^{\deg(M)}
          \end{eqnarray}
      where, if we order the elements of $G(\Jl)$
            lexicographically as $\Ms$, then for $i=1, \ldots, m$
          $$\set(M_i)=\{j \in \{1, \ldots, n\} \st x_j \in (M_1,
            \ldots, M_{i-1}):M_i \}.$$

          In our case, as the degree of each of the monomials
          generating $\Jl$ is $b+1$, Equation (\ref{e:poincare}) turns
          into
           \begin{eqnarray}\label{e:poincare2}
             \begin{array}{ll}
              P_{R/\Jl}(q,t)& \displaystyle =1 +\sum_{M
               \in G(\Jl)}(1+qt)^{|\set(M)|}qt^{b+1}\\
              &\\
              & \displaystyle =1 + qt^{b+1}\sum_{M
                \in G(\Jl)}(1+qt)^{|\set(M)|}.
              \end{array}
           \end{eqnarray}

       So now we focus on $|\set(M_j)|$ for $M_j \in
           \{\Ms\}$. Suppose $M_j=x_{i_1}\cdots x_{i_{b+1}}$, where
           $i_1 < \ldots < i_{b+1}$. Then each $M \in \{M_1, \ldots,
           M_{j-1}\}$ is of the form $M=x_{u_1}\cdots x_{u_{b+1}}$,
           with $u_1 < \ldots < u_{b+1}$, and $M$ is lexicographically
           smaller than $M_j$. So the relationship between the indices
           is such that
           \begin{center} $u_1 < i_1$ \hspace{.5in} or
                \hspace{.5in} if $u_1=i_1, \ldots, u_l=i_l$
                 then $u_{l+1} < i_{l+1}$.
           \end{center}

        So by an argument identical to that in the proof of Lemma~\ref{l:I-lemma}
        $$\begin{array}{ll}
        \set(M_j)&=\{u \leq n \st x_u \in (M_1,
            \ldots, M_{j-1}):M_j \}\\
         &\\
         &=\{u \leq i_{b+1} \st x_u \ndiv M_j\}.
      \end{array}$$
         We can now conclude that \begin{equation}\label{e:setM}
         |\set(M_j)|= i_{b+1}-(b+1). \end{equation}

      We have shown that, if $M$ is any degree $b+1$ square-free
      monomial with highest index $u$ (that is, $x_u \ddiv M$ and
      $x_v\ndiv M$ for $v >u$), then $|\set(M)|=u-(b+1)$. So to
      compute the sum in~(\ref{e:poincare2}), all we have to
      do is count the number of square-free degree $b+1$ monomials
      with highest index $u$, for any given $u$. This number is
      clearly ${u-1}\choose{b}$. So for a given $i$, the number of
      degree $b+1$ square-free monomials $M$ with $|\set(M)|=i$ is
      exactly ${b+i}\choose{b}$.

        Therefore
        \begin{eqnarray}\label{e:poincare3}
         \begin{array}{ll}
          P_{R/\Jl}(q,t)& = \displaystyle 1 + qt^{b+1}
            \sum_{i=0}^{n-b-1}{{b+i} \choose{b}}(1+qt)^i\\
           &\\
           &\displaystyle = 1 + qt^{b+1}\sum_{i=0}^{a}{{b+i}
            \choose{b}}(1+qt)^i
          \end{array}
         \end{eqnarray}
    since by  Equation (\ref{e:setM}), $i$ can reach at most $n-b-1$, which by definition is equal to $a$.

        \item[{\em Step 2.}]  Since $\El$ is generated by a regular
             sequence over $R/\Jl$ (Proposition~\ref{p:e-regular}), we
             can use a mapping cone construction to find its minimal
             graded resolution (see
             Discussion~\ref{d:mapping-cones}). We do this by adding
             the generators of $\El$, one at a time, to $\Jl$, and
             applying Corollary~\ref{c:mapping-cone}. As the
             generators $e_1(\xs), \ldots, e_b(\xs)$ of $\El$ have degrees
             $1,\ldots, b$, respectively, each time we add a
             $e_i(\xs)$, the Poincar\'e series gets multiplied by a factor
             of $(1+qt^i)$, and hence

              $$\begin{array}{lll}
               P_{R/\Il}(q,t)& = \displaystyle \prod_{k=1}^{b}(1+qt^k)\cdot
                P_{R/\Jl}(q,t)&\\

              &\displaystyle = \prod_{k=1}^{b}(1+qt^k) \cdot\, \Big( 1
               + q t^{b+1} \sum_{i=0}^a { b+i \choose b} (1+q\,t)^{i}
               \Big)& \ ({\rm from~(\ref{e:poincare3})}).

          \end{array}$$

     \end{itemize}

          \end{proof}

\begin{corollary}[The set of generators of $\Il$ is minimal] Let 
$\mu=(a \ba b)$ be a hook. The generating set for $\Il$ described in
Proposition \ref{p:generators} is minimal.
\end{corollary}

   \begin{proof} The number of generators of $\Il$ is by definition 
     ${n\choose{b+1} }+ b$. On the other hand, the minimal number of
     generators of $\Il$ is the first Betti number $\be_1$ of $R/\Il$,
     which is the coefficient of $q$ in the Poincar\'e series
     $P_{R/\Il}(q,1)$.  It is easy to see by Theorem~\ref{mainthm}
     that this coefficient is $$ b+1 +\sum_{i=1}^a {b+i \choose b}.$$
     So all we have to show is that $${n\choose{b+1} }+ b= b+1
     +\sum_{i=1}^a {b+i \choose b}$$ which is equivalent to showing
     that $${n\choose{b+1} }=\sum_{i=0}^{n-b-1} {b+i \choose b}.$$
     This last equation follows easily from induction on $n$.
   \end{proof}

\subsection{Combinatorial interpretations}\label{s:recursive}

From Theorem \ref{mainthm} it can be seen that the Poincar\'e series
of $\Il$ can be defined recursively.  For simplicity, for any hook
$\mu=(a\mid b)$ we denote by $P_{(a\mid b)}(q,t)$ the bigraded
Poincar\'e series $P_{R/\mathcal{I}_{(a\mid b)}}(q,t)$, and by
$P_{(a\mid b)}(q)$ the nongraded Poincar\'e series
$P_{R/\mathcal{I}_{(a\mid b)}}(q,1)$.

We start with the vertical partition $(0\mid b)$. In this case the ideal $\Il$ is generated only by the elementary symmetric functions; the quotient of $\Il$ is the coinvariant algebra, a well-known representation of the symmetric group (see e.g.,  \cite{Hu}).  The nongraded Poincar\'e series in this case is

\begin{align*}
    &  \begin{matrix}
               \Young{
                \Box{ \circ }\cr
                 \Box{ \circ }\cr
                \Box{ \circ }\cr
                \Box{}  \cr }
     \end{matrix}
    && P_{(0\mid b)}(q) = (1+q)^{b+1}.
\end{align*}

Using Equation (\ref{e:mainthm}), by subtracting  $P_{(a-1\mid b)}(q)$ from  $P_{(a\mid b)}(q)$, we find that the nongraded Poincar\'e polynomial  of $\mu$ satisfies the following recurrence:

$$P_{(a\mid b)}(q) = P_{(a-1\mid b)}(q)+ {a+b \choose b} q (1+q)^{a+b}.$$

This recursion allows us to compute the nongraded Poincar\'e polynomial of $(a\mid b)$ by adding one cell at a time to the first row of the vertical partition $(0\mid b)$

{\allowdisplaybreaks
\begin{align*}
   &  \begin{matrix}
               \Young{
                \Box{ \circ }\cr
                 \Box{ \circ }\cr
                \Box{ \circ }\cr
                \Box{}  & \Box{\bullet} \cr }
     \end{matrix}
     && P_{(1\mid b)}(q) = (1+q)^{b+1} \Big( 1 + { b+1 \choose b } q \Big)\\
      &  \begin{matrix}
               \Young{
                \Box{ \circ }\cr
                 \Box{ \circ }\cr
                \Box{ \circ }\cr
                \Box{}  & \Box{\bullet}   & \Box{\bullet} \cr }
     \end{matrix}
     && P_{(2\mid b)}(q) = (1+q)^{b+1} \Big( 1 + { b+1 \choose b } q  + { b+2 \choose b} q (1+q) \Big),
\end{align*}
}

until we reach the hook $(a\mid b)$, which gives us
$$P_{(a\mid b)}(q) = (1+q)^{b+1} \Big( 1 + q \sum_{i=1}^a { b+i \choose b} (1+q)^{i-1} \Big).$$

The graded Poincar\'e polynomial satisfies a similar recurrence:
\begin{align*}
P_{(0\mid b)}(q,t)&= \prod_{k=1}^{b+1} (1+qt^k) \\ P_{(a\mid b)}(q,t)&
 =P_{(a-1\mid b)}(q,t) + \prod_{k=1}^{b} (1+qt^k) \,\cdot\, q t^{b+1}
 {b+a \choose a} (1+qt)^a.
\end{align*}
Once again, like the nongraded case, one can use this recurrence to build $P_{(a\mid b)}(q,t)$ starting from $P_{(0\mid b)}(q,t)$.

 Now we turn to the question of a combinatorial interpretation of the coefficients $\beta_{i,j}$ of $P_{(a\mid b)}(q,t)$.
In the case of the vertical partition $(0\mid b)$ such an interpretation is given by Cauchy's $t$-binomial theorem, which states that
$$\prod_{k=1}^n (1+q\, t^k) = \sum_{k=0}^n q^k \, t^{\frac{k(k+1)}{2}} {
n \choose k} _t.$$
Here
 \[{ n \choose k} _t:= \frac{[n]_t!}{[k]_t! \, [n-k]_t!} \] are the $t$-binomial coefficients which have many interesting combinatorial interpretations~(\cite{St1}), and
\begin{eqnarray}\label{e:symbol}
[j]_t !:=[1]_t \, [2]_t \, \cdots
[j]_t {\rm \; \; with \;\; } [j]_t: = 1+t+\ldots +t^{j-1} .\end{eqnarray}

The following question begs to be answered: in general, is it
possible to find a combinatorial interpretation for the graded Betti
numbers of the De Concini-Procesi ideals?

\idiot{The following was a formula suggested by Francois, but it is not correct as stated.

(Generating series for $P_{(a\mid b)}(q,t)$)
$$\sum_{a,b\geq 0} P_{(a\mid b)}(q,t) z^a w^b= \frac{(1+q)(1-z-w(1+q))}{(1-z)(1-t(1+q))(1-(z+w)(1+q))}.$$
}

\section{Regularity of Hooks}\label{s:regularity}

\begin{definition}[Castelnuovo-Mumford regularity]\label{d:regularity}
Let $I$ be an ideal of a $R=k[\xs]$. The \emph{Castelnuovo-Mumford
regularity} or simply \emph{regularity} of $R/I$, denoted by
$\reg(R/I)$ is defined as the maximum value of of $j-i$ where the
graded Betti number $\beta_{i,j}\neq 0$ in a minimal free resolution
of $R/I$.
\end{definition}

\begin{corollary}[Regularity of hooks]
Let $\mu=(a \mid b)$ be a hook. Then $\reg(R/I)=b(b+1)/2$.
\end{corollary}

    \begin{proof}  The graded Betti numbers $\beta_{i,j}$ appear as the
     coefficients of the Poincar\'e series
     $$P_{R/\Il}(q,t)=\underbrace{\prod_{k=1}^{b}
      (1+qt^k)}_{\mbox{Factor 1}} \,\cdot\, \underbrace{\Big( 1 + q
      t^{b+1} \sum_{i=0}^a { b+i \choose b} (1+q\,t)^{i}
      \Big)}_{\mbox{Factor 2}}.$$ So the question is to find the term
      $q^it^j$ in this polynomial, where the coefficient $\beta_{i,j}$
      is nonzero and $j-i$ is maximum.  The terms with nonzero
      coefficients in each factor are of the following forms:
      $$\begin{array}{lll}
       \mbox{Factor 1:} & q^m t^{b_1+\ldots+b_m} & \mbox{where }
                1\leq b_1<\ldots<  b_m\leq b, \ 0\leq m \leq b, \\
       \mbox{Factor 2:} & q^{e+1} t^{e+b+1} & \mbox{where }0\leq e \leq a.
      \end{array}$$

       To show that $\displaystyle \reg(R/I)=\frac{b(b+1)}{2}$, we need
       to show that this bound is achieved by the possible choices of
       $j-i$, and is the maximum possible bound.

       Consider the terms in Factor 1. We have
       $$\begin{array}{ll}
        \displaystyle b_1+\ldots+b_m -m &\leq \Big((b-(m-1))+(b-(m-2))+\ldots + b\Big) -m\\
        &\\
        &=\Big((1+2+\ldots+b)-(1+2+\ldots+(b-m))\Big) -m \\
        &\\
        &\displaystyle =\frac{b(b+1)}{2} - \frac{(b-m)(b-m+1)}{2} -m \\
        &\\
        &\displaystyle \leq \frac{b(b+1)}{2}.
         \end{array}$$

       Similarly, for terms in Factor 2, since $b\geq 0$, we have
       $$e+b+1-(e+1)=b \leq \frac{b(b+1)}{2}.$$

       For the product of a term in Factor 1 and a term in Factor 2 we have
       $$\begin{array}{ll}
         b_1+\ldots+b_m+e+b+1 -(m + e+1) & \\
        &\\
        = b_1+\ldots+b_m+b-m &\\
        &\\
        \displaystyle \leq \frac{b(b+1)}{2} - \frac{(b-m)(b-m+1)}{2}+b-m &
         \mbox{ same argument as above}\\
         &\\
        \displaystyle = \frac{b(b+1)}{2}- (b-m)\left(\frac{b-m+1}{2}
        -1\right) &\\
        &\\
        \displaystyle = \frac{b(b+1)}{2}-(b-m)\left(\frac{b-m-1}{2}\right) &\\
        &\\
        \displaystyle \leq \frac{b(b+1)}{2} & \mbox{\ \ since\ \ } m\leq b.
     \end{array}$$
      The bound is achieved if $m=b$, so that $b_1=1, \ldots,b_m=b$, and
      for any $e$, so that we have the term with nonzero coefficient
     $$q^{e+b+1}t^{(1+\ldots+b)+e+b+1}=q^{e+b+1}t^{\frac{b(b+1)}{2}+e+b+1}$$
     which clearly has the property that $j-i=\displaystyle
     \frac{b(b+1)}{2}$, as desired.
    \end{proof}

\section{The Hilbert series of hooks}\label{s:hilbert}
The goal of this section is to find the Hilbert
 series of $R/\Il$ when $\mu$ is a hook partition, namely, the series
 $$h_{R/\Il}(q)=\sum_{s=0}^{\infty}\dimn_k(R/\Il)_sq^s,$$
where as usual $\dimn_k$ means dimension as a vector space over $k$.
This has been done in the general case of a partition $\mu=(\mu_1,\ldots,\mu_n)$ of $n$ by  Garsia and Procesi. In~\cite{GP}, they provide an explicit basis for $R/\Il$ as a $\mathbb{Q}$-module, from which it follows that
\begin{equation}
{\rm dim}_{\mathbb{Q}}(R/\Il)={n \choose \mu_1,\ldots,\mu_n}.  \label{dimension}
\end{equation}
and
\begin{equation}\label{hilbert-series-GP}
h_{R/\Il}(q)=\sum_{\lambda \vdash n} f^\lambda \; K_{\lambda \mu}(1/q) \;q^{n(\mu)}.
\end{equation}
Here, $f^\lambda$ and $n(\mu)$ are two well-known parameters associated with partitions (\cite{St2}), and $K_{\lambda \mu}(q)$ are the Kostka-Foulkes polynomials we referred to in the introduction (\cite{LS}). The computation of $K_{\lambda \mu}(q)$ is somewhat complicated. This motivates us to use the results of this paper to give a new description of the Hilbert series in the case of hooks.

Let $\mu=(a \ba b)$ be a hook partition of $n$, and consider the
 ideal $\Il=\Jl+\El$. Since $R/\Jl$ is a Cohen-Macaulay ring (Corollary~\ref{c:J-CM}), and
the generators $e_1(\xs), \ldots, e_b(\xs)$ of $\El$ form a regular sequence over
$R/\Jl$ (Proposition~\ref{p:e-regular}), it follows that (see~\cite{V}
Theorem~4.2.5)
\begin{equation}\label{e:hilbertseries}
h_{R/\Il}(q)=\prod_{i=1}^b (1-q^i)h_{R/\Jl}(q).
\end{equation}

\idiot{(Sara) Also, see~\cite{V} Theorem~2.2.7. The actual theorem that
states this is due to Stanley (\cite{V} refers to this). The version
that appears in Eisenbud's book (pp. 551 or so) is wrong-- A and B
should be exchanged.}

So we focus on finding $h_{R/\Jl}(q)$. Recall that $\Jl$ is generated by all square-free monomials of degree
$b+1$ with variables in $\{\xs\}$. So each graded piece $(R/\Jl)_s$ is generated
by all monomials of degree $s$, involving $c$ of the $n$ variables with $c \leq b$.  There are $n\choose c$ ways of choosing $c$ variables from $\{\xs\}$. Choose such a monomial, without loss of generality,
$$x_1^{a_1} \ldots \, x_c^{a_c}.$$ We need to choose the positive
integers $a_1, \ldots, a_c$ such that $a_1+ \ldots + a_c=s$.

This is classically equivalent to inserting $c-1$ bars between
the sequence of integers $1, \ldots, s$, as below:
$$ 1,\ldots, a_1 \ba a_1+1, \ldots, a_1+a_2 \ba \ldots \ba
(a_1+\ldots+a_{c-1})+1, \ldots, s.$$

What we are doing here is choosing $c-1$ of the $s-1$ available slots,
and there are ${s-1}\choose{c-1}$ ways of doing that. So we have

$$h_{R/\Jl}(q) = 1+ \sum_{s=1}^{\infty} \sum_{c=1}^{b}
{n \choose c}
{s-1 \choose c-1}
q^s.$$

Therefore, by Equation (\ref{e:hilbertseries})
$$\begin{array}{ll}
h_{R/\Il}(q) & = \displaystyle \prod_{i=1}^{b}(1-q^i) \left ( 1+ \sum_{s=1}^{\infty} \sum_{c=1}^{b} {n \choose c} {s-1 \choose c-1} q^s\right  )\\
&\\
& = \displaystyle \prod_{i=1}^{b}(1-q^i) \left ( 1+ \sum_{c=1}^{b} {n \choose c} \sum_{s=1}^{\infty} {s-1 \choose c-1} q^s\right  ).
  \end{array}
$$
On the other hand
$$\begin{array}{lll}
\displaystyle \sum_{s=1}^{\infty} {s-1 \choose c-1} q^s & \displaystyle =\sum_{s=c}^{\infty} {s-1 \choose c-1} q^s&\\
&&\\
&= \displaystyle q^c\sum_{s=c}^{\infty} {(s-c)+(c-1) \choose c-1} q^{s-c} &\\
&&\\
&= \displaystyle q^c\sum_{s=c}^{\infty} {(s-c)+(c-1) \choose s-c} q^{s-c} &
\displaystyle {\rm \  because \ }{i+j \choose i}={i+j \choose j}\\
&&\\
&= \displaystyle \frac{q^c}{(1-q)^c}  & {\rm \  because \ } \displaystyle \sum_{j=0}^{\infty}{i+j \choose j}q^j=\frac{1}{(1-q)^{i+1}}\ \  (\cite{W}).
\end{array}
$$

So
 $$\begin{array}{ll}
h_{R/\Il}(q) & \displaystyle = \prod_{i=1}^{b}(1-q^i)\left ( 1+ \sum_{c=1}^{b} {n \choose c}\frac{q^c}{(1-q)^c} \right ) \\
&\\
&\displaystyle = \frac{(1-q)^b}{(1-q)^b}\prod_{i=1}^{b}(1-q^i)\left ( 1+ \sum_{c=1}^{b} {n \choose c}\frac{q^c}{(1-q)^c} \right ) \\
&\\
&\displaystyle = \prod_{i=1}^{b}\left ( \frac{1-q^i}{1-q} \right )
\sum_{c=0}^{b} {n \choose c}q^c(1-q)^{b-c}  \\
&\\
&\displaystyle = [b]_q!\  \sum_{c=0}^{b} {n \choose c}q^c(1-q)^{b-c}  \\
  \end{array}
$$
where $[b]_q!$ is described in (\ref{e:symbol}) above.
So we have proved that
\begin{proposition}
Let $\mu=(a\mid b)$ be a hook partition of $n$. Then
\begin{equation}\label{final-formula}
h_{R/\Il}(q)=[b]_q!\  \sum_{c=0}^{b} {n \choose c}q^c(1-q)^{b-c}.
\end{equation}
\end{proposition}
Note that if we set $q=1$ in (\ref{final-formula}), we find that
$$\dimn_k (R/\Il)=\frac{n!}{(a+1) !} =\frac{n!}{\mu_1 !}$$
as was expected by Formula (\ref{dimension}).

\todo{Also (see, for example,~\cite{BH})
$$h_{R/\Il}(q)=\frac{P_{R/\Il}(-1,q)}{(1-q)^n}$$ where
$P_{R/\Il}(q,t)$ is the Poincar\'e series given in our paper (this is
messy at the moment)!

In \cite[Corollary 1.2]{E2}, there is a specific formula given for the
Hilbert function using the graded Betti numbers. Here is what it says
($m$= length of resolution, $r+1=$number of variables, and each free
module in the resolution $F_i=\bigoplus_jR(-a_{i,j})$:
$$H_M(d)=\sum_{i=0}^m(-1)^i\sum_j {{r+d-a_{i,j}}\choose{r}}.$$

This is probably what we looked at before (in fact the same as above),
but check again to be sure.  See also Cor. 1.10.}


\end{document}